\documentclass[a4paper,12pt]{amsart}
\usepackage{amssymb}
\usepackage{amsthm, amsmath, latexsym, amsfonts}
\usepackage[all]{xy}
\newtheorem{Lemma}{Lemma}
\newtheorem{Proposition}{Proposition}
\newtheorem{Theorem}{Theorem}

\newcommand{\Pro}{\mathrm{Pr}}
\newcommand{\Bl}{\mathrm{Bl}}

\newcommand{\N}{\mathbb{N}}

\newcommand{\Q}{\mathbb{Q}}
\newcommand{\C}{\mathbb{C}}
\newcommand{\bP}{\mathbb{P}}

\title[Cohomology of moduli spaces of higher genus stable maps]{Towards the cohomology of moduli spaces \\ of higher genus stable maps}

\author{Claudio Fontanari}

\email{claudio.fontanari@polito.it}\curraddr{
Dipartimento di Matematica \\ Politecnico di Torino \\
Corso Duca degli Abruzzi 24 \\ 10129 Torino \\ Italy.}

\begin{document}

\begin{abstract}
We prove that the orbifold desingularization of the moduli space 
of stable maps of genus $g = 1$ recently constructed by Vakil and 
Zinger has vanishing rational cohomology groups in odd degree 
$k < 10$.
\end{abstract}

\maketitle

\section{}
Let $\overline{\mathcal{M}}_{g,n}(\bP^r,d)$ be the moduli space of stable maps 
of degree $d$ from $n$-pointed curves of genus $g$ to the projective space 
$\bP^r$. For a comprehensive introduction to this fashinating geometric 
object we refer to the classical text \cite{FP:97}, which not only 
contains a careful construction of Kontsevich moduli spaces but also outlines 
their crucial application to Gromov-Witten theory. 

For $g=0$ it turns out that $\overline{\mathcal{M}}_{0,n}(\bP^r,d)$ 
carries a natural structure of smooth orbifold, in particular its 
rational cohomology is well-defined and well-behaved. Indeed, in 
the last few years there has been a flurry of research about the 
cohomological properties of $\overline{\mathcal{M}}_{0,n}(\bP^r,d)$
after the first steps moved in \cite{BF:02}: it is worth mentioning 
at least the contributions \cite{BO:03}, \cite{BF:03}, \cite{Oprea:04},
\cite{Oprea:05}, \cite{Oprea:06}, \cite{MM:04},\cite{MM:05}, \cite{MM:06}, 
\cite{MM:06bis}, \cite{GeP:05}, \cite{Cox:05}, \cite{Cox:05bis}, 
\cite{Cox:05ter}.

For higher genus $g > 0$, instead, $\overline{\mathcal{M}}_{g,n}(\bP^r,d)$
can be arbitrarily singular (see \cite{Vakil:06}) and it is in general 
nonreduced with several components of exceptional dimension. As a consequence, 
a cohomological investigation of the underlying topological spaces has never 
been addressed and it seems to be completely out of reach. However, at least in 
the case $g=1$, the beautiful construction performed in \cite{VZ:06} has 
recently opened a new frontier to the research in the field. 

Namely, in \cite{VZ:06} it is shown that the closure 
$\overline{\mathcal{M}}^0_{1,n}(\bP^r,d)$ of the stratum 
$\mathcal{M}^0_{1,n}(\bP^r,d)$ corresponding to stable 
maps with smooth domain allows an orbifold desingularization 
$$
\widetilde{\mathcal{M}}^0_{1,n}(\bP^r,d) \to \overline{\mathcal{M}}^0_{1,n}(\bP^r,d)
$$
which can be explicitely described as a sequence of blow-ups along 
smooth centers. Moreover, the natural $(\C^*)^{r+1}$-action on 
$\overline{\mathcal{M}}^0_{1,n}(\bP^r,d)$ lifts to 
$\widetilde{\mathcal{M}}^0_{1,n}(\bP^r,d)$ and the 
corresponding fixed loci are explicitely determined. 
As a by-product, the rational cohomology groups 
$$H^k(\widetilde{\mathcal{M}}^0_{1,n}(\bP^r,d), \Q)$$ 
are well-defined and can be approached via localization techniques. 

In the present paper we take this direction and prove the following result: 

\begin{Theorem}\label{odd}
We have $H^k(\widetilde{\mathcal{M}}^0_{1,n}(\bP^r,d), \Q) = 0$ 
for every $n \ge 0$, $r \ge 1$ and odd $k < 10$.
\end{Theorem}

We point out that the natural surjective map 
$$
\widetilde{\mathcal{M}}^0_{1,n}(\bP^r,d) \to \overline{\mathcal{M}}_{1,n}
$$ 
together with the well-known fact that
$$
H^{11}(\overline{\mathcal{M}}_{1,11}, \Q) \ne 0
$$
(see, for instance, \cite{Pikaart:95}, proof of Corollary~4.7, or \cite{GrP:01}, Proposition~2) show that in the above statement the bound on the degree $k$ is 
sharp. 

Our proof is based on a suitable version of the Bialynicki-Birula decomposition 
for singular varieties (\cite{Kirwan:88} and \cite{Oprea:04}), which reduces the 
claim of Theorem~\ref{odd} to the analogous claim on the fixed locus. 
Next, a careful analysis of its connected components allows us to conclude 
thanks to previous results on the cohomology of the moduli spaces of 
stable curves of genus $g \le 1$ (\cite{Keel:92} and \cite{Getzler:98}). 

We work over the field $\C$ of complex numbers. 

Many thanks are due to Gilberto Bini and Barbara Fantechi for stimulating 
conversations about this research project. 

\section{}

We begin by recalling some more or less well-known facts about 
the Bialynicki-Birula decomposition for orbifolds.

\begin{Proposition}\label{BB}
Let $X$ be a smooth projective orbifold with a $\C^*$-action, let $F$ 
be the fixed locus and let $F_i$ denote its connected components. 
Then

(i) $X$ is the disjoint union of locally closed subvarieties $S_i$ 
such that every $S_i$ retracts onto the corresponding $F_i$ and 
$$
\overline{S}_i \subseteq \bigcup_{j \ge i} S_j
$$ 

(ii) the Betti numbers of $X$ are 
$$
h^m(X) = \sum_i h^{m-2n_i}(F_i)
$$
where $n_i$ is the codimension of $S_i$. 
\end{Proposition}

\proof (i) follows from \cite{Kirwan:88}, (1.14) and (1.15), pp. 388-389; 
(ii) follows from (i) and \cite{Oprea:04}, Lemma~6~(i). 

\qed

We also need the following technical result: please refer to \cite{VZ:06}, 
Section~2.1, for all non standard definitions and notation. 

\begin{Lemma}\label{blow}
Fix $N \in \N$. Let $Y$ be a smooth (orbi)variety and let $\{X_1, \ldots, X_n \}$ 
be a properly intersecting collection of subvarieties of $Y$ such that 
$\bigcap_{i \in I} X_i$ is smooth for every $I \subseteq \{1, \ldots, n \}$
and $H^k(\bigcap_{i \in I} X_i, \Q)$ $=0$ for every odd $k < N$.
Let $Y_1 :=Y$ and for $i=1, \ldots, n$ let $f_i: Y_{i+1} \to Y_i$ be the 
blow-up of the proper transform of $X_i$ in $Y_i$. 

(i) If $X_i^j$ denotes the proper transform of $X_i$ in $Y_j$, then 
$\bigcap_{i \in I} X_i^j$ is smooth for every $j \ge 1$ and
every $I \subseteq \{j \le i \le n \}$ and 
$H^k(\bigcap_{i \in I} X_i^j, \Q)$ $=0$ for every odd $k < N$. 

(ii) If moreover $H^k(Y, \Q)=0$ for every odd $k < N$, then $H^k(Y_i, \Q)=0$ 
for every odd $k < N$ and for every $i \le n+1$.  
\end{Lemma}

\proof (i) By induction on $j$, the case $j=1$ holding by assumption. We have
$$
\bigcap_{i \in I} X_i^j = \bigcap_{i \in I} \Pro_{X_{j-1}^{j-1}} X_i^{j-1} = 
\Pro_{X_{j-1}^{j-1}} \bigcap_{i \in I} X_i^{j-1} = \Bl_Z X
$$
with $Z := \bigcap_{i \in I} X_i^{j-1} \cap X_{j-1}^{j-1}$ and 
$X := \bigcap_{i \in I} X_i^{j-1}$. Here the first equality holds 
by definition, the second one by \cite{VZ:06}, Lemma 2.3~(2), and 
the third one by \cite{VZ:06}, Lemma 2.3~(1). 

Now, it is well-know how to compare cohomology groups after a blow-up 
along a smooth subvariety of an orbifold (see \cite{GH:78}, p.~605 
and Proposition on p.~606, and \cite{KL:04}, footnote on p.~514, 
for the orbifold case).   

In particular, it follows that if both $H^k(X, \Q)=0$ and $H^k(Z, \Q)=0$ 
for odd $k < N$, then $H^k(\Bl_Z X, \Q)=0$ for odd $k < N$. Hence 
in our case we obtain by induction both the smoothness and the vanishing 
result. 

(ii) By induction on $i$, the case $i=1$ holding by assumption. We have 
$$
Y_{i+1} = \Bl_{X_i^i} Y_i 
$$
with $H^k(Y_i, \Q) =0$ for odd $k < N$ by induction and $H^k(X_i^i, \Q) =0$ 
for odd $k < N$ by (i). Hence the claim follows from the cohomological 
properties of blow-ups recalled above. 

\qed

\emph{Proof of Theorem~\ref{odd}.} By \cite{VZ:06}, Theorem~1.1,  
$\widetilde{\mathcal{M}}^0_{1,n}(\bP^r,d)$ is a smooth projective 
orbifold and the natural $(\C^*)^{r+1}$-action on the projective 
space $\bP^r$ lifts to an action on  
$\widetilde{\mathcal{M}}^0_{1,n}(\bP^r,d)$. 
By Propostion~\ref{BB}~(ii), in order to prove Theorem~\ref{odd} 
it is enough to show that the odd cohomology of the fixed loci 
vanishes in degree $k < 10$. 

By \cite{VZ:06}, Section~1.4, such fixed loci are either 
\begin{equation}\label{first}
\prod \overline{\mathcal{M}}_{g,m}
\end{equation}
with $g \le 1$ and $m \in \N$ (see \cite{VZ:06}, p.~15, and \cite{MirSym}, p.~541), or
\begin{equation}\label{second}
\prod \overline{\mathcal{M}}_{0,i} \times \bP^j \times 
\widetilde{\mathcal{M}}_{1, (I,J)}
\end{equation}
with $i, j \in \N$ (see \cite{VZ:06}, p.~19), where $\widetilde{\mathcal{M}}_{1, (I,J)}$ 
is constructed in \cite{VZ:06}, pp.~26--27. Roughly speaking, 
$\widetilde{\mathcal{M}}_{1, (I,J)}$ is obtained from the moduli 
space $\overline{\mathcal{M}}_{1, I \sqcup J}$ 
of stable curves of genus $1$ by successively 
blowing up certain subvarieties $\overline{\mathcal{M}}_{1, \rho}$ 
(where $\rho$ is a suitable index) 
and their proper transforms (see \cite{VZ:06}, p.~25). 
We just point out a crucial fact: from the definition of 
$\overline{\mathcal{M}}_{1, \rho}$ (see \cite{VZ:06}, p.~23)
it follows that 
$$
\overline{\mathcal{M}}_{1, \rho} = \overline{\mathcal{M}}_{1, h} \times 
\prod \overline{\mathcal{M}}_{0, k}
$$
with $h, k \in \N$. More generally, if 
$\bigcap_{\rho_\alpha} \overline{\mathcal{M}}_{1, \rho_\alpha} \ne \emptyset$
for some subset of indices $\{ \rho_\alpha \}$ then 
\begin{equation}\label{third}
\bigcap_{\rho_\alpha} \overline{\mathcal{M}}_{1, \rho_\alpha} = 
\overline{\mathcal{M}}_{1, s} \times \prod \overline{\mathcal{M}}_{0, t}
\end{equation}
with $s, t \in \N$. 

Now, recall that all odd degree cohomology of $\overline{\mathcal{M}}_{0, n}$ 
vanishes by Keel's results (see \cite{Keel:92}) and that 
$H^k(\overline{\mathcal{M}}_{1, n}, \Q) = 0$ for odd $k < 10$ 
according to Getzler (see \cite{Getzler:98} and \cite{GrP:01}, eq.~(4)). 
Hence the required vanishing trivially follows for fixed loci of type
(\ref{first}) and is reduced to the case of 
$\widetilde{\mathcal{M}}_{1, (I,J)}$ for fixed loci of type 
(\ref{second}). On the other hand, by (\ref{third}) we can 
apply Lemma~\ref{blow}~(ii) to our situation with $N = 10$, 
$Y = \overline{\mathcal{M}}_{1, I \sqcup J}$, 
$X_i = \overline{\mathcal{M}}_{1, \rho}$, and 
$Y_{n+1} = \widetilde{\mathcal{M}}_{1, (I,J)}$, 
so once again Getzler's claims allow us to conclude
and the proof is over.

\qed

\end{document}